\documentclass{amsart}

\usepackage{amsmath,amsfonts,amssymb,amsthm,amscd,latexsym,euscript, stmaryrd}
\usepackage[all]{xy}

\newcommand{\op}{{\ensuremath{\textup{op}}}}
\newcommand{\Hor}{\ensuremath{\textup{Hor}}}

\newcommand{\dgrm}[1]{\ensuremath{\smash{\underset{\widetilde{\hphantom{#1}}}{#1}} \mathstrut}}

\newcommand{\cal}[1]{\ensuremath{\mathcal #1}}

\newtheorem {theorem1}{Theorem}[section]
\newtheorem {theorem}[theorem1]{Theorem}

\newtheorem {proposition}[theorem1]{Proposition}
\newtheorem {lemma}[theorem1]{Lemma}
\theoremstyle{definition}
\newtheorem {definition}[theorem1]{Definition}
\newtheorem {example}[theorem1]{Example}
\theoremstyle{remark}
\newtheorem {remark}[theorem1]{Remark}

\newcommand{\cat}[1]{\ensuremath{\EuScript #1}}

\newcommand{\Sets}{\ensuremath{\cal{S}}\emph{ets} }

\newcommand{\colim}{\ensuremath{\mathop{\textup{colim}}}}

\def\Lan{\ensuremath{\textup{Lan}}}

\newcommand{\rarrow}{\rightarrow}

\newcommand{\Rosicky}{Rosick\'y }

\newcommand{\simpsets}{\cal S}
\newcommand{\CC}{\cal C}
\newcommand{\DD}{\cal D}
\newcommand{\sS}{\cal S}
\newcommand{\sm}{_{\text{sm}}}
\newcommand{\smallFun}{\sS^{\DD}\sm}
\newtheorem*{VariableNoNum}{{\VariableText}}
\newenvironment{titled}[1]
     {\def\VariableText{{#1}}\begin{VariableNoNum}}
     {\end{VariableNoNum}}
\newcommand{\LS}{\cal L}
\newcommand{\SetFor}[1]{\cal O_{#1}}
\newcommand{\MapOf}[1]{{\cal Map}(#1)}
\newcommand{\MapsSDD}{\MapOf{\sSDD}}
\newcommand{\sSsSop}{\sS^{\sS^\op}}
\newcommand{\sSsS}{\sS^{\sS}}
\newcommand{\MapsS}{\MapOf \sS}
\newcommand{\MapsSeq}{\MapOf\sS_{\text{eq}}}
\newcommand{\sSDD}{\sS^{\DD}}
\newcommand{\Glue}{\operatorname{\Gamma}}
\newcommand{\glue}{\operatorname{\gamma}}
\newcommand{\Imap}{\operatorname{Inc}}
\begin{document}

\SelectTips{cm}{10}

\title  [Homotopy theory of small diagrams]
        {Homotopy theory of small diagrams over large categories}

\author{Boris Chorny}
\author{William G.~Dwyer}
\thanks{
The second author was partially supported by National Science
Foundation grant DMS-0204169. The authors would like to thank
the Mittag-Leffler Institute for its hospitality during the
period in which this paper was completed. 
}
\address{D-MATH, ETH Zentrum, 8059 Zurich, Switzerland}
\address{Department of Mathematics, University of Notre Dame, Notre Dame, IN 46556, USA}

\email{chorny@math.ethz.ch}
\email{dwyer.1@nd.edu}



\date{\today}
\dedicatory{}
\commby{}

\begin{abstract}
  Let $\cal D$ be a large category which is cocomplete. We construct a
  model structure (in the sense of Quillen) on the category of small
  functors from $\cal D$ to simplicial sets.  As an application we
  construct homotopy localization functors on the category of
  simplicial sets which satisfy a stronger universal property than the
  customary homotopy localization functors do.
\end{abstract}

\maketitle

\section{Introduction}
Let $\simpsets$ be the category of simplicial sets.  In this paper we
introduce axiomatic homotopy theory into the study of functors from a
large category $\cal D$ into $\simpsets$, in other words, into the
study of diagrams in $\simpsets$ indexed by $\cal D$.  Such diagrams
arise naturally (for instance in the treatment of Goodwillie calculus,
which in one form deals with functors from $\simpsets$ itself into
$\simpsets$) but in the past they have been dealt with by \emph{ad
  hoc} techniques.  The novelty of of our approach is the introduction
of a model category structure, which allows for the use of standard
tools from axiomatic homotopy theory.

There is an obvious set--theoretic difficulty in dealing with the
index categories we wish to consider: if $\cal D$ is large, the
totality of natural transformations between two functors $\sS\to\sS$
does not necessarily form a set, and so the collection of all such
functors is not even a category in the usual sense, much less a model
category. We overcome this difficulty by restricting our attention to
the category $\smallFun$ of \emph{small} (\ref{DefineSmall}) functors
$\DD\to\sS$.  This category is always cocomplete. If $\DD$ itself is
cocomplete, then $\smallFun$ is also complete, and it is in this
situation we can construct a model category on $\smallFun$. This model
structure reduces to the ordinary projective model category structure
on the category of all functors $\DD\to\sS$ if $\DD$ is small; observe
though that for technical reasons our model structure in general lacks
functorial factorization.

We discuss in detail two examples, $\DD=\sS^\op$ and $\DD=\sS$.  For
$\DD=\sS^\op$, we generalize the arguments of \cite{DK} to show that
our model structure on $\smallFun$ is Quillen equivalent to the
equivariant model structure developed by Farjoun \cite{Farjoun} on the
category of maps of spaces.  The model category $\sS^{\sS}\sm$ does
not seem to have an analogous interpretation. The category of
pro-spaces may be viewed as dual to the subcategory of
pro-representable functors in $\sS^{\sS}\sm$, and its model structure
\cite{EH76} \cite{Isaksen-strict} is perhaps the closest relative to
our model structure on $\sS^{\sS}\sm$.

An immediate application of the model structure on $\sS^{\sS}\sm$ is a
construction of homotopy localizations in this category. Although this
construction itself involves factorizations and is thus non-functorial
in $\sS^{\sS}\sm$, an application of the constuction to the identity
functor yields an object of $\sS^{\sS}\sm$ (i.e., a functor
$\sS\to\sS$) which is equivalent to the ordinary homotopy localization
functor on simplicial sets but has a stronger universal property.  We
finish the paper by using these homotopy localization functors to
construct natural $A$-Postnikov towers in~$\sS$. Another application of the
model structure on $\sS^{\sS}\sm$ is developed in \cite{BCR}.

\subsection*{Acknowledgements}
The first author would like to thank Jir\'i \Rosicky for helpful
conversations on the early stages of this project. We are deeply
obliged to Steve Lack for sharing with us his unpublished work.

\subsection{Notation} We continue to let $\sS$ denote the category of
simplicial sets, which we also refer to as the category of
\emph{spaces}. If $\CC$ and $\DD$ are categories, we simplify notation
by using $\CC^{\DD}$ to denote the category $\CC^{\DD}\sm$ of all
small functors $\DD\to\CC$. If $\DD$ itself is small, this is the
category of all functors $\DD\to\CC$. A \emph{simplicial category} is
a category enriched over~$\sS$, such as $\sS$ itself; functors between
two such categories are assumed to respect the enrichments, in the
sense that they provide simplicial maps between the respective
function complexes.

\section{Preliminaries on small functors}\label{prelim}

The object of study of this paper is homotopy theory of functors from
a large simplicial category to~$\sS$. The totality of these functors
does not form a category in the usual sense, since the natural
transformations between two functors need not form a set in general,
but rather a proper class. We are willing to be satisfied with a
treatment of a reasonable subcollection of functors, a subcollection
which does form a category. The purpose of this section is to describe such a
subcollection.

\begin{definition}\label{DefineSmall}\label{DefineRepresentable}
  Let $\DD$ be a (not necessarily small) simplicial category. A
  functor $\dgrm X:\DD\to\sS$ is \emph{representable} if there is an object
  $D\in \DD$ such that $\dgrm X$ is naturally equivalent to $R^D$, where
  $R^D(D')=\hom_{\DD}(D, D')$. A functor $\dgrm X\colon \cat C \rarrow
  \cal S$ is called \emph{small} if \dgrm X is a small weighted
  colimit of representables.
\end{definition}


\begin{remark}\label{SmallIsKanExtension}
  Since the category of small functors is tensored over simplicial
  sets, the \emph{small weighted colimit} above may be expressed as a coend of the form 
  \begin{equation}\label{DisplayWeightedColimit}
    R^F\otimes_{\cat I}G = \int^{I\in\cat I}R^{F(I)}\otimes G(I), 
  \end{equation}
  where $\cat I$ is a small category and $F\colon \cat I\to \DD$,
  $G\colon \cat I \to\sS$ are functors. Here $R^F\colon \cat I^\op\to
  \cal S^\DD$ assigns to $I\in\cat I$ the representable functor
  $R^{F(I)}\colon \DD\to\sS$. For the general treatment of weighted
  limits and colimits see \cite{Kelly}.

  Since the simplicial tensor structure on the category of small
  fuctors $\sS^\DD$ is given by the objectwise direct product, we
  will use $\dgrm X\times K$ to denote tensor product of $\dgrm X\in
  \cal S$ with $K\in \sS$.

  The above coend is the (enriched) left Kan extension of the functor
  $G$ over the functor $F$. Using the transitivity of left Kan
  extensions, it is easy to see that the following four conditions are
  equivalent \cite[Prop.~4.83]{Kelly}:
  \begin{itemize}
  \item $\dgrm X\colon\DD\to\sS$ is a small functor,
  \item there is a small simplicial category $\cat I$ and a functor $G:\cat
    I\to \sS$, such that $\dgrm X$ is isomorphic to the left Kan
    extension of $G$ over some functor $\cat I\to\DD$,
  \item there a small simplicial subcategory $i\colon\DD'\to\DD$ 
    and a functor $G:\DD'\to\sS$,
    such that $\dgrm X$ is isomorphic to the left Kan extension of 
    $G$ over $i$, and
  \item there is a small full simplicial subcategory $i\colon\DD_{\dgrm
      X}\to
    \DD$ such that $\dgrm X$ is isomorphic to the left Kan extension
    of $i^*(\dgrm X)$ over~$i$.
  \end{itemize}
\end{remark}

  If $D\in \DD$ and $\dgrm Y$ is a functor $\DD\to\sS$, then by
  Yoneda's lemma the simplicial class of natural transformations
  $R^D\to\dgrm Y$ is $\dgrm Y(D)$; in particular, this simplicial
  class is a simplicial set. It follows easily that if $\dgrm X$ is a
  small functor $\DD\to\sS$, then the natural tramsformations $\dgrm
  X\to\dgrm Y$ also form a simplicial set (this also follows from
  \ref{SmallIsKanExtension} above and the adjointness property of the
  left Kan extension). In particular, the
  collection of all small functors is a simplicial category. 

\begin{remark}\label{SmallFunctorCategory}
  M.G.~Kelly \cite{Kelly} calls small functors \emph{accessible} and
  weighted colimits \emph{indexed}. He proves that small functors
  form a simplicial category which is closed under small (weighted)
  colimits \cite[Prop.~5.34]{Kelly}. 
\end{remark}

In order to do homotopy theory we need to work in a category which is
not only cocomplete, but also complete (at least under finite limits).
Fortunately, there is a simple sufficient condition in the situation
of small functors. 

\begin{theorem}
  If $\DD$ is cocomplete, then the category $\sS^{\DD}$ of small
  functors $\DD\to\sS$ is complete.
\end{theorem}

\begin{remark}
  There is a long story behind this theorem.  P.~Freyd \cite{Freyd}
  introduced the notion of \emph{petty} and \emph{lucid} set-valued
  functors. A set-valued functor is called \emph{petty} if it is a
  quotient of a small sum of representable functors. Any small functor
  is clearly petty. A functor $F\colon \cat A \rarrow \Sets$ is called
  \emph{lucid} if it is petty and for any functor $G\colon \cat
  A\rarrow \Sets$ and any pair of natural transformations $\alpha,
  \beta\colon G\rightrightarrows F$, the equalizer of $\alpha$ and
  $\beta$ is petty.  Freyd proved \cite[1.12]{Freyd} that the category
  of lucid functors from $\cat A^\op$ to \Sets is complete if and only
  if $\cat A$ is \emph{approximately} complete (that means that the
  category of cones over any small diagram in \cat A has a weakly
  initial set).  J.~\Rosicky then proved \cite[Lemma 1]{Rosicky} that
  if the category $\cat A$ is approximately complete, a functor
  $F\colon \cat A^\op \rarrow \Sets$ is small if and only if it is
  lucid.  Finally, these results were partially generalized by B.~Day
  and S.~Lack \cite{Lack} to the enriched setting. They show, in
  particular, that the category of small $\cat V$--enriched functors
  $\cat K^\op\to \cat V$ is complete if \cat K is complete and \cat V is a
  symmetric monoidal closed category which is locally finitely
  presentable as a closed category.  This last condition is certainly
  satisfied if $\cat V=\sS$.
\end{remark}

\section{A model category on $\cal S^{\DD}$}\label{MainSection}
As usual, $\sS^{\DD}$ denotes
the category of small functors $\DD\to\sS$.

\begin{theorem}\label{MainTheorem}
Assume that $\DD$ is cocomplete. Then the category  $\cal S^{\DD}$
has a model category structure in which weak equivalences and
fibrations are defined objectwise and the cofibrations are the maps
which have the left lifting property with respect to acyclic
fibrations. (The factorizations provided by this model category
structure are not necessarily functorial.)
\end{theorem}

\begin{remark}
The use of ``objectwise'' above signifies that a maps $F\to G$ is a
weak equivalence (fibration) if and only if for each $X\in\DD$ the
induced map $F(X)\to G(X)$ is a weak equivalence (fibration) of
simplicial sets. We are using the ordinary model category structure on
simplicial sets, in which a map is a weak equivalence if its geometric
realization is a weak equivalence of topological spaces,  and a
fibration if it is a Kan fibration (see, e.g., \cite[Thm.~3.6.5]{Hovey}).
\end{remark}

Recall from \ref{DefineRepresentable} the notion of representable
functor, as well as the notation $R^D=\hom(D,\text{--})$ for the
functor represented by~$D$. We first need a definition and some
lemmas, which exhibit yet additional uses of word \emph{small}.

\begin{definition}
  A collection $\LS$ of objects in a category $\cat B$ is said to be
  \emph{locally small} in $\cat B$ if for every object $X$ of $\cat B$
  there exists a \emph{set} of objects $\SetFor X\subset\LS$ such that
  any map $Y\to X$ with $Y\in \LS$ can be factored as a composite
  $Y\to Y'\to X$ for some~$Y'\in\SetFor X$.
\end{definition}

\begin{remark}\label{CoSolution}
  More standardly, the statement that $\LS$ is locally small is expressed by
  saying that $\LS$ satisfies the \emph{co--solution set} condition.
  The set $\SetFor X$ is called the co--solution set associated
  to~$X$. Our terminology follows \cite{Farjoun}, since the idea of
  the proof of Theorem~\ref{MainTheorem} also goes back to \cite{Farjoun}.
\end{remark}

\begin{lemma}\label{loc-small}
The collection of representable functors is locally small 
in $\sS^{\DD}$.
\end{lemma}

\begin{proof}
Suppose that \dgrm X is in $\sS^{\DD}$, and write \dgrm X as a small
weighted colimit as in (\ref{DisplayWeightedColimit}).
Given a representable functor $R^D$, consider the simplicial set
$\hom(R^D, \dgrm X)$. From the generalized Yoneda lemma, and the fact
that weighted colimits of diagrams are computed levelwise, we obtain:
\begin{equation}
\hom(R^D, \dgrm X) = \dgrm X(D)= R^F(D)\otimes_{\cat I} G 
= \int^{I\in\cat I}\hom(R^D, R^{F(I)})\times G(I).
\end{equation}
Comparing the sets of the zero simplices of the simplicial sets above,
we conclude that every map $R^D\rarrow \dgrm X$ factors through a map
$R^D\rarrow R^{F(I)}$ for some $I\in\cat I$.
\end{proof}

\begin{definition}\label{DefineMapCat}
The \emph{category of maps in $\sSDD$}, denoted $\MapsSDD$, is the category
whose objects are the arrows $f\colon \dgrm X\to\dgrm Y$ in $\sSDD$. A morphism
$f\to f'$ is a commutative diagram
\[
\xymatrix{
\dgrm X \ar[d]^f \ar[r]     &       \dgrm X' \ar[d]^{f'}\\
\dgrm Y \ar[r]                     &       \dgrm Y'
}
\]
\end{definition}

\begin{lemma}\label{MapSmall}
  If $g:K\to L$ is a map of spaces, then the collection
  \[\LS(g)=\{R^D\times K\to R^D\times L \mid D\in\DD\}\] is locally small in
  $\MapsSDD$.
\end{lemma}

\begin{remark}\label{UnionSmall}
  It follows immediately that if $\{g_\alpha\}_{\alpha\in A}$ is a set
  of maps between spaces, then the union $\cup_\alpha\LS(g_\alpha)$ is
  also locally small in $\MapsSDD$.
\end{remark}

\begin{titled}{Proof of \ref{MapSmall}}
Consider a morphism 
\[
\xymatrix{
R^D\times K \ar[d]^{g_D} \ar[r]     &       \dgrm X \ar[d]^f\\
R^D\times L \ar[r]                     &       \dgrm Y
}
\]
in $\MapsSDD$. This morphism gives rise, by adjunction, to the
following commutative diagram:
\[
\xymatrix{
R^D \ar@/^/[drr] \ar@/_/[ddr] \ar@{-->}[dr]^\varphi&&\\
    & \dgrm W \ar[d] \ar[r]     &       \dgrm X^K \ar[d]\\
    & \dgrm Y^L \ar[r]          &       \dgrm Y^K,
}
\]
where $\dgrm W$ is defined so that the square is a pullback square. By
Lemma \ref{loc-small} there exists a set of representable functors
$\SetFor {\dgrm W}$, such that any morphism from a representable
functor to \dgrm W can be factored through an object in $\SetFor{\dgrm
  W}$. Now take $\SetFor f = \{F \times K\to F\times L \mid
F\in\SetFor {\dgrm W}\}$, and observe that any map from $g_D$ to $f$
will factor, by adjunction, through one of the objects in $\SetFor f$.
\qed
\end{titled}

Let us now briefly remind the setup of the generalized small object
argument, which applies for locally small collection of maps with
small domains. The reader might want to consult \cite{pro-spaces} 
for a more extensive discussion. Suppose that $\LS$ is a locally
small collection in $\MapsSDD$, that $f\colon\dgrm X\to \dgrm Y$ is an
object in $\MapsSDD$, and that $\SetFor f$ is the associated
co--solution set for~$f$ (\ref{CoSolution}). We define
$\Glue^1_{\LS}(f)$ to be the natural map $\glue^1_{\LS}(f)\to Y$, where
$\glue^1_{\LS}(f)$ is determined by the following pushout diagram:
\[\xymatrix{
\coprod_\beta \dgrm U_\beta \ar[d] \ar[r]  & \dgrm X\ar[d]\\
\coprod_\beta \dgrm V_\beta \ar[r] & \glue^1_{\LS}(f). 
}
\]
Here $\beta$ runs through the set of pairs $(g_\beta,h_\beta)$, where
$g_\beta\colon \dgrm U_\beta\to \dgrm V_\beta$ belongs to $\SetFor f$ and
$h_\beta\colon g_\beta\to f$ is a morphism in $\MapsSDD$. It is easy to see that
the map $\dgrm X\to \dgrm Y$ extends to a map $\glue^1_{\LS}(f)\to
\dgrm Y$. For $n>1$, we let $\glue^n_{\LS}(f)=\glue^1_{\LS}(\Glue^{n-1}_{\LS}(f))$,
and $\Glue^n_{\LS}(f)\colon\glue^n_{\LS}(f)\to Y$ the induced natural map.
Finally, $\glue^\infty_{\LS}(f)$ denotes $\colim_n\glue^n_{\LS}(f)$, and
$\Glue^\infty_{\LS}(f)\colon \glue^\infty_{\LS}(f)\to Y$ is the evident natural map.

Recall that a simplicial set is said to be \emph{finite} if it has a
finite number of nondegenerate simplicies and a finite simplicial set
$K$ is $\aleph_0$--small in the category of simplicial sets, in the sense that
$\hom(K,\text{--})$ commutes with countable sequential colimits.

In order to conclude, by the generalized small object argument, that
the induced map $\Glue^\infty_{\LS}(f):\glue^\infty(f)\to \dgrm Y$ has
the right lifting property with respect to all of the maps in $\LS$,
the class $\LS$ must satisfy an additional condition (to local
smallness) that all domains of maps in $\LS$ are $\lambda$--small for
some fixed cardinal $\lambda$.

Yoneda's lemma and smallness of finite simplicial sets imply the last
condition if $\LS$ is the collection $\cup_\alpha\LS(g_\alpha)$ for a
set $\{g_\alpha\colon K_\alpha\to L_\alpha\}_{\alpha\in A}$ of
monomorphisms between finite simplicial sets. $\LS$ is locally small
by \ref{UnionSmall}.

\begin{remark}
  The construction of $\Glue^1_{\LS}(f)$ or $\Glue^\infty_{\LS}(f)$ is
  not functorial 
  unless the co--solution sets
  $\SetFor f$ depend in some natural way on~$f$. This would be the
  case, for instance, if  $\SetFor f=\LS$ for all~$f$, but of course this
  would be allowed only if $\LS$ itself is a set. Another example
  where $\SetFor f$ depends functorially on $f$ occurs in the
  equivariant model category of \cite{Farjoun}. See \cite{PhDI} for
  the construction of functorial factorizations in this model
  category. In general, there are two versions of the generalized
  small object argument: functorial and non-functorial
  \cite{pro-spaces}. We apply the non-functorial version in this work. 
\end{remark}

\begin{proof}[Proof of \ref{MainTheorem}]
There are several variations in the literature of the definition of a model category. We prove the axioms MC0--MC5 in the form of \cite{DS}.
  The required limits and colimits exist by the discussion in the
  previous section. The `2-out-of-3' axiom and the fact that weak
  equivalences are closed under retracts follow from the corresponding
  properties of the category \cal S. By the definition of cofibration,
  every cofibration has the left lifting property with respect to any
  acyclic fibration.

  Although our model category is not cofibrantly generated
  \cite[Section~2.1]{Hovey}, it has a similar structure, namely, it is
  class-cofibrantly generated \cite[Def.~1.3]{pro-spaces}. In order to verify the second lifting property and two factorization
  properties, let us define the classes of generating cofibrations and
  generating acyclic cofibrations to be
\begin{align*}
I &= \{R^D\times \partial \Delta^n \hookrightarrow R^D\times \Delta^n | D\in \cal D, n\geq 0\} \\
J &= \{R^D\times \Lambda^n_k \hookrightarrow R^D\times \Delta^n | D\in
\cal D, n>0, n\geq k\geq 0\}
\end{align*}
where as usual $\Delta^n$ is the $n$-simplex, $\partial\Delta^n$ its
boundary, and $\Lambda^n_k$ the space obtained by removing the $k$'th
face of $\Delta^n$ from $\partial\Delta^n$.  A map in $\sS$
is an acyclic fibration if and only if it has the right lifting
property with respect to the inclusions $\partial
\Delta^n\to\Delta^n$, $n\ge0$ and so it follows by adjunction that a
map in $\sSDD$ is an acyclic fibration if and only if it has the right
lifting property with respect to the maps in~$I$. Similarly, a map in
$\sSDD$ is a fibration if and only if it has the right lifting property
with respect to the maps in~$J$.

Suppose that $f\colon \dgrm X\to\dgrm Y$ is a map in $\sSDD$, and note that by
\ref{UnionSmall} discussion above the classes $I$ and $J$ permit the
generalized small object argument \cite{pro-spaces}. Therefore, the
composite $\dgrm X\to \glue^\infty_I(f)\to \dgrm Y$ is a factorization
of $f$ into the composite of cofibration with an acyclic fibration,
while $\dgrm X\to \glue^\infty_J(f)\to \dgrm Y$ is a factorization
into the composite of an acyclic cofibration and a fibration. 

The second lifting property is achieved by a standard trick; see,
e.g., \cite{DK}. Given an acyclic cofibration $f:\dgrm A \to \dgrm B$,
let $\dgrm C=\glue^\infty_J(f)$ and factor $f$ as a composite $\dgrm
A\to \dgrm C\to\dgrm B$. By construction the map $\dgrm A\to \dgrm C$ has
the left lifting property with respect to any fibration. Since $\dgrm
C \to\dgrm B$ is actually an acyclic fibration (by the `2-out-of-3'
property), lifting in the diagram
\[
\xymatrix{
         \dgrm A \ar[d]\ar[r] & \dgrm C\ar[d]\\
          \dgrm B \ar[r]      &  \dgrm B
}
\]
shows that $\dgrm A \to \dgrm B$ is a retract of $\dgrm A \to \dgrm C$
and thus also has the left lifting property with respect to any
fibration.
\end{proof}

\section{Example: $\DD  = \cal S^\op$}
In order to illustrate the concept of the projective model structure
on the category of small functors $\cal S^{\DD}$ by a familiar model
category, we consider $\DD= \sS^\op$. In this case we construct a
Quillen equivalence between $\cal S^{\cal D}$ and the category
$\MapsSeq$; the subscript ``eq'' signifies that this is the category
$\MapsS$ of maps in $\sS$ (\ref{DefineMapCat}) endowed with the
``equivariant'' or ``fine'' model structure of \cite{Farjoun}.

Recall from \cite{DZ} that an \emph{orbit} in $\MapsS$ is a
diagram $A\to B$ in $\sS$ whose colimit is isomorphic to the
one--point space $*=\Delta^0$. Since the colimit of $A\to B$ is $B$,
an orbit in $\MapsS$ is simply an object of the form $A\to *$, and so, via
the functor which assigns to such a diagram the space $A$, the
category $\cal O$ of orbits is equivalent to $\sS$.  We will let
$O_A=(A\to *)$ be the orbit in $\MapsS$ corresponding the
space~$A$. The following definition was given in \cite{Farjoun}:

\begin{definition}\label{EquivModelStr}
  The \emph{equivariant model structure} or \emph{fine model
  structure} on $\MapsS$ is the model category $\MapsSeq$ in which the
  underlying category is $\MapsS$, and in which a map $X\to Y$ between
  objects of $\MapsS$ is a weak equivalence (fibration) if and only
  for each $A\in \sS$  it gives a weak equivalence (fibration)
  $\hom(O_A, X)\to \hom(O_A,Y)$ in $\sS$.
\end{definition}


The above definition suggests assigning to each object $X$ of $\MapsS$ the
diagram $X^{\cal O}\colon \sS^\op\to\sS$ sending $A$ to $\hom(O_A,X)$; the
functor $(\text{--})^{\cal O}\colon \MapsS \to \sS^{\sS^\op}$ both preserves
and reflects weak equivalences. 

\begin{lemma}
  If $X\in\MapsS$, the $X^{\cal O}\colon \sS^\op\to\sS$ is a small functor;
  in particular, $X\mapsto X^{\cal O}$ gives a functor $\MapsS\to\sSsSop$.
\end{lemma}

\begin{proof}
  This was essentially shown in \cite{DF}. More precisely, Farjoun
  proved (in a more general context) that for any object $X\in\MapsS$
  there exists a small full subcategory $i\colon \cal O_{ X}
  \hookrightarrow \cal O$ such that $ X^{\cal O}$ is the left Kan
  extension of $i^\ast(X^{\cal O})$ (cf. \ref{SmallIsKanExtension}).
\end{proof}

Let us construct the left adjoint to this functor by verifying the
conditions of adjoint functor theorem: the orbit-point functor
obviously preserves limits, so we have to verify the solution-set
condition. This means for any small diagram $\dgrm X\in \cal S^{\cal
S^\op}$ we need to find a set of arrows $f_i\colon\dgrm X \rarrow
Y_i^{\cal O}$ such that any arrow $f\colon \dgrm X\rarrow
Z^\cal O$, for $Z\in \MapsS$ factorizes as $f=(k)^{\cal O}\circ
f_i$. 

Recall from \cite{DK} that for every full simplicial subcatetgory
$\cat I \subset \cal O \cong \sS$ there is a pair of ajoint functors
\begin{equation}\label{D-K}
|\text{--}|_{\cat I}\colon \cal S^{\cat I^\op} \leftrightarrows
\MapsS :\! (\text{--})^{\cat I},
\end{equation}
(which is a Quillen equivalence if one considers the model
structure induced by the set \cat I of orbits on $\MapsS$).

If \dgrm X is small, then it is a left Kan extension of  $i^*\dgrm X\colon
\cat I^\op\rarrow \cal S$ for a small simplicial full subcategory \cat I
of the orbit category $\cal O\cong \sS$ and $i\colon \cat I^\op\to \cal O^\op$.

Given $f\colon \dgrm X\to Z^{\cal O}$, consider $i^*f\colon
i^*\dgrm X \to i^* Z^{\cal O} = Z^\cat I$ and look at the
adjoint map $(i^*f)^\sharp\colon |i^*\dgrm X|_\cat I\to Z$. 
Let $Y = |i^*\dgrm X|_\cat I$ and $k = (i^*f)^\sharp$. Consider
the unit of adjunction (\ref{D-K}): $i^*\dgrm X\to (|i^*\dgrm X|_\cat I)^\cat I =
Y^\cat I = i^*(Y^{\cal O})$ and apply the left Kan extension
along $i$. Composing with the counit of the adjunction between the
left Kan extension and $i^*$ we obtain the map $f'\colon \Lan_i
i^*\dgrm X = \dgrm X \to \Lan_i i^*(Y^{\cal O}) \to 
Y^{\cal O}$ and the required factorization $\dgrm X\rarrow
Y^{\cal O}\rarrow Z^{\cal O}$, so that $f= k^{\cal O}f'$.  

That means that there exists the left adjoint to the functor $(\text{--})^{\cal
O}$, we will call it \emph{realization} functor and denote the realization of \dgrm Z by $|\dgrm Z|$.

More explicitly, the left adjoint to the functor $X\mapsto X^{\cal O}$ is
the functor which assigns to $\dgrm Y\in\sSDD$ the coend
$\Imap\times_{\DD}\,\dgrm Y$, where $\Imap:\DD^\op=\sS\to\sS$
is the identity inclusion (for notational reasons, let
$\DD=\sS^\op$). Of course on the face of it this is a large coend, but
it actually gives a functor for the following reason. Since $\dgrm Y$
is a small diagram, there is a small full subcategory $i:\cat 
I\subset\DD$ such that $\dgrm Y$ is the left Kan extension of
$i^*(\dgrm Y)$ over $i$. It then follows from associativity properties of
coends that $\Imap\times_{\DD}\,\dgrm Y$ is isomorphic to the small
coend $\Imap\times_{\cat I}\,\,i^*(\dgrm Y)$.


\begin{proposition}
The functors $X\mapsto X^{\cal O}$ and $\dgrm Y\mapsto |\dgrm Y|$ form
a Quillen pair, which give a Quillen equivalence between $\sSsSop$
and $\MapsSeq$.
\end{proposition}


\begin{proof}
  To produce the Quillen pair, it will suffice by
  \cite[8.5.3]{Hirschhorn} to show that any (acyclic) fibration in
  $\MapsSeq$ is preserved by the functor $X\mapsto X^{\cal O}$; this,
  however, follows immediately from definition of the model category
  structures on $\MapsSeq$ (\ref{EquivModelStr}) and on $\sSsSop$
  (\ref{MainTheorem}).

  In order to show that the pair of functors is a Quillen equivalence
  we have to prove that for any cofibrant diagram $\dgrm Y \in \cal
  S^{\cal S^\op}$ and for any fibrant $X\in \MapsS$, a
  map $f\colon\dgrm Y\rarrow  X^{\cal O}$ is a weak equivalence
  if and only if the adjoint map $f^\sharp\colon |\dgrm Y|\rarrow
  X$ is a weak equivalence. But, by definition of weak
  equivalences, $f^\sharp$ is a weak equivalence if and only if
  $(f^\sharp)^{\cal O}\colon |\dgrm Y|^{\cal O}\rarrow  X^{\cal
  O}$ is a weak equivalence, so it will suffice, by the 2-out-of-3
  property, to show that the unit of the adjunction induces a weak
  equivalence $g\colon \dgrm Y\rarrow |\dgrm Y|^{\cal O}$ for every
  cofibrant diagram \dgrm Y. This we now do.

  From the small object argument (\S\ref{MainSection}) we know that \dgrm Y is a retract
  of $\Glue^\infty_I(\emptyset\to\dgrm Y)$, where $\emptyset$ is the
  empty diagram.  Since retracts are preserved by all functors, and
  retracts of weak equivalences are weak equivalences, we can assume
  that $\dgrm Y=\Glue^\infty_I(\emptyset\to\dgrm Y)$. Let $\dgrm
  Y_i=\Glue^n_I(\emptyset\to\dgrm Y)$; then
\[
\dgrm Y = \colim(\dgrm Y_1\rarrow\dgrm Y_2\rarrow\cdots\rarrow \dgrm Y_n\rarrow\cdots),
\]
where $\dgrm Y_n\rarrow \dgrm Y_{n+1}$ is obtained by a pushout
\[
\xymatrix{
\coprod_\alpha R^{A_\alpha}\otimes \partial\Delta^n \ar[d] \ar[r]     &       \dgrm Y_n \ar[d]\\
\coprod_\alpha R^{A_\alpha}\otimes \Delta^n \ar[r]                     &       \dgrm Y_{n+1}\,.
}
\]
Since left adjoints commute with colimits we obtain:
\[
|\dgrm Y| = \colim(|\dgrm Y_1|\rarrow|\dgrm Y_2|\rarrow\cdots\rarrow |\dgrm Y_n|\rarrow\cdots),
\]
where $|\dgrm Y_n|\rarrow |\dgrm Y_{n+1}|$ is obtained by a pushout
\[
\xymatrix{
\coprod_\alpha |R^{A_\alpha}|\otimes \partial\Delta^n \ar[d] \ar[r]     &       |\dgrm Y_n| \ar[d]\\
\coprod |R^{A_\alpha}|\otimes \Delta^n \ar[r]                     &       |\dgrm Y_{n+1}|.
}
\]
But $|R^{A}| \cong O_A=(A\to *)$, as can be verified by using the
above coend description of the realization functor, or by a simple
adjointness verification.
Hence, the pushout diagram above becomes
\[
\xymatrix{
\coprod_\alpha O_{A_\alpha}\otimes \partial\Delta^n \ar[d] \ar[r]     &       |\dgrm Y_n| \ar[d]\\
\coprod_\alpha O_{A_\alpha} \otimes \Delta^n \ar[r]                     &       |\dgrm Y_{n+1}|.
}
\]

But it was shown in \cite{Farjoun} and \cite{PhDI} that the functor
$(-)^{\cal O}$ commutes up to a weak equivalence with all colimits
of the above form. This immediately leads the conclusion that the
natural map $\dgrm Y\rarrow |\dgrm Y|^{\cal O}$ is a weak equivalence
in the $\cal S^{\cal S^\op}$.
\end{proof}

\section{Application: homotopy localization of spaces}
In this section we take $\cal D = \cal S$ and we present an
application of the projective model structure on 
$\cal S^{\cal S}$ to homotopy localization in the category of
spaces.

We first recall some basic notions. 
Suppose that $f\colon A\rarrow B$ is a cofibration of spaces. A space
$Z$ is said to be \emph{$f$-local} if $Z$ is fibrant and
$f^*:\hom(B,Z)\to\hom(A,Z)$ is a weak equivalence in $\sS$; a map $X\to Y$ is
an \emph{$f$-equivalence} if $\hom(Y,Z)\to \hom(X,Z)$ is a weak
equivalence in $\sS$ for every $f$-local~$Z$. Finally, an
$f$-equivalence $X\to X'$
is an \emph{$f$-localization map} if $X'$ is $f$-local.

It is well known (see \cite{Bous:factor} and \cite{Farjoun-book}) that
there exists a homotopy idempotent, coaugmented, simplicial homotopy functor
$L_f\colon \cal S \rarrow \cal S$ which has the following two
properties:
\begin{enumerate}
\item for any $X\in \cal S$ the
coaugmentation $\eta_X\colon X\to L_f(X)$ is an $f$-localization map, and
\item for
every map $g\colon X \rarrow Z$, where $Z\in \cal S$ is $f$-local,
there exists a factorization of $g$
\[
\xymatrix{
X \ar[dr]_{\eta_X} \ar[rr]^g   &     &  Z\\
            & L_f X\ar@{-->}[ur]_h,
}
\]
and in this factorization the map $h$ is unique up to simplicial homotopy.
\end{enumerate}

We produce a localization functor which is weakly equivalent to the
one above, but which has a stronger universal property. Assume as
usual that $f\colon A\to B$ is a cofibration of spaces.

\begin{theorem}\label{functorial_localization}
  There here exists a homotopy idempotent, coaugmented, small,
  simplicial, homotopy functor $L_f:\sS\to\sS$ with the following two
  properties 
\begin{enumerate}
  \item for any $X\in \cal S$ the coaugmentation $\eta_X\colon X\to
  L_f(X)$ is an $f$-localization map, and
  \item for every coaugmented  functor $\dgrm L\colon\sS\to \sS$
  taking $f$-local values, there exists a factorization 
\[
\xymatrix{
\mathrm{Id} \ar[dr]_{\eta} \ar[rr]^\varepsilon   &     &  \dgrm L.\\
            & L_f \ar@{-->}[ur]_\zeta
}
\]
In this factorization the natural transformation $\zeta$ is unique
up to a simplicial homotopy (of natural transformation).
  \end{enumerate}
\end{theorem} 

\begin{proof}
  We sketch the proof, with references.  Given a cofibration $f\colon
  A\hookrightarrow B$, consider the class $N$ of maps in $\cal S^{\cal
    S}$ given by $N = \{ R^C\times f \;|\;  C\in \cal S\}$.
  Then $N$ is locally small in $\MapOf\sSsS$ (\ref{MapSmall}), and
  just as in the fixed-point-wise situation of \cite{PhDI}, the class
\begin{align*}
  \Hor(N) = \left\{\left.(R^C\times f) \boxempty \left( \underset{\Delta^n} {\overset
  {\partial\Delta^n} {\downarrow}}\right) \,\right| \, C\in \sS, n\geq
  0\right\} \\ =
  \left\{\left. R^C \times A\times \Delta^n\coprod_{R^C \times A\times\partial\Delta^n}
  R^C \times B\times\partial\Delta^n\to R^C \times B\times \Delta^n \,\right|\, C\in \sS, n\geq 0 \right\}.
\end{align*}
  permits the generalized small object argument
  \cite{pro-spaces}. (Note that although the functors
  $\hom(A,\text{--})$ and $\hom(B,\text{--})$ do not necessarily commute with sequential
  colimits in~$\sS$, they do commute with well--ordered colimits of
  sufficiently high cofinality.) Observe now that the identity functor
  $\mathrm{Id}$ and the constant functor $*$ on $\sS$ are small and in
  fact representable; one is $R^C$ for $C=*$ and the other for $C$
  being the empty diagram. Therefore, taking $\LS = \Hor(N)\cup J$ and
  applying the generalized small object technique, we can factor the
  map $ \mathrm{Id} \rarrow *$ into a composite $\mathrm{Id}\rarrow
  \dgrm K \rarrow *$ in which, by construction, $\dgrm K$ is a small
  simplicial functor. 

  There are three properties of this factorization to notice.
  First, by the very nature of the generalized small object technique,
  the map $\dgrm K\to *$ has the right lifting
  property with respect to the maps in $\LS$. This amounts to the
  assertion that for each $C\in\sS$ the space $\dgrm K(C)$ is fibrant
  and has the right lifting property with respect to the maps in
  $\Hor(N)$, i.e., to the assertion that $\dgrm K(C)$ is $f$-local. 
  Secondly, if $\dgrm L\colon \sS\to\sS$ is a functor which takes on
  $f$-local values, then the induced map $\hom(\dgrm K,
  \dgrm L)\to\hom(\mathrm{Id}, \dgrm L)$ is an acyclic fibration in~$\sS$. This
  follows from the way in which $\dgrm K$ is constructed from iterated
  pushouts of the maps in $\LS$, and the fact that for any $g\colon
  U\to V$ in $\Hor(N)$ and any $f$-local space $Z$, the restriction
  map $\hom(V,Z)\to\hom(U,Z)$ is an acyclic fibration. By picking
  $C\in \sS$ and applying this observation to the coextended diagram $\dgrm L$
  given by $\dgrm L(X)=\hom(\hom(X,C), Z)$, one sees that for any
  $f$-local space $Z$, $\hom(K(C),Z)\to\hom(C, Z)$ is an acyclic
  fibration. In particular $C\to K(C)$ is an $f$-equivalence.
  Finally, $\dgrm K$ is a homotopy functor; in fact the above
  considerations show that $\dgrm K$ is an $f$-localization functor,
  and it follows formally from the definition that such functors take
  weak equivalences between spaces into simplicial homotopy
  equivalences between fibrant spaces. To finish the proof, it is
  enough to take $L_f=\dgrm K$. We leave it to the reader to deduce from
  \ref{functorial_localization}(2) that $L_f$ is homotopy idempotent.
\end{proof}

\section{Example: A functorial $A$-Postnikov tower}
It is well known that there exists a construction of the classical
Postnikov tower, which is functorial `as a tower'. However, this
construction, due to Moore, is \emph{ad hoc} and does not allow for a
natural generalization. Our new construction of localizations 
provides a general method of obtaining functorial towers. 

\begin{example}
E.~Farjoun discussed a Postnikov tower with respect to a space $A$
\cite{Farjoun-book}. This is a construction that associates to every
space $X$ a tower of spaces $\cdots\to P_{\Sigma^n A}X\to
P_{\Sigma^{n-1} A}X\to \cdots \to P_{\Sigma A}X\to P_{A}X$, where $P_B
= L_{B\to\ast}$ is the nullification functor. The classical
construction of localizations ensured that each level in this tower is
a functor in $X$, but not the whole tower. We take an advantage of
localizations with the functorial universal property in order to
obtain an equivalent tower functorial in $X$. 

Let $f_n$ be the map $\Sigma^n A\to \ast$ for all $n\geq 0$. From now
on denote by $P_{\Sigma^n A} = L_{f_n}$ the localization functor with
functorial universal property constructed in the previous section. A
fibrant simplicial set $X$ is $f_n$-local iff $\ast \simeq
\hom(\Sigma^n A, X) = \hom(\Sigma^{n-1} A, \Omega X) =
\Omega\hom(\Sigma^{n-1} A, X)$. Therefore if a fibrant $X$ is
$f_{n-1}$-local, then $X$ is $f_n$-local. By
Theorem~\ref{functorial_localization} for each $n>0$ there exists a
natural transformation $\zeta_n\colon L_{f_n}\to
L_{f_{n-1}}$. Combining $\zeta_n$ for all $n>0$ we obtain the required
tower of functors $\cdots\to P_{\Sigma^n A}\overset
{\zeta_n}\longrightarrow P_{\Sigma^{n-1} A}\overset
{\zeta_{n-1}}\longrightarrow \cdots \to P_{\Sigma A}\overset
{\zeta_1}\longrightarrow P_{A}$. If $A = S^0$, then we recover a new
construction of the classical Postnikov tower. 
\end{example}

\begin{remark}
The construction of the localization functor with a stronger
functorial property can be carried out also in the stable model
cateogry of spectra. As an application we can obtain the functorial
decomposition of spectra into a tower of chromatic layers.
\end{remark}

\bibliographystyle{abbrv}
\bibliography{Xbib}

\end{document}